\documentclass[letterpaper,11pt]{article}

\usepackage[utf8]{inputenc}
\usepackage{authblk}
\usepackage{amssymb}
\usepackage{amsmath}
\usepackage{amsxtra}
\usepackage{amsfonts}
\usepackage{mathtools}
\usepackage{hyperref}
\usepackage{bm}
\usepackage{siunitx}
\usepackage{tikz}
\usepackage{fullpage}
\usepackage{tkz-euclide}
\usepackage{listings}
\usepackage{algorithm}
\usepackage{graphicx} 
\usepackage{booktabs}
\usepackage{fancyhdr}
\usepackage{algorithm}
\usepackage{algpseudocode}
\usepackage{enumitem}
\usepackage{amsthm}

\usepackage{graphicx}
\usepackage{caption}
\usepackage{subcaption}

\usepackage{amssymb}
\usepackage{amsmath}
\usepackage{amsthm}
\usepackage{graphicx}
\usepackage{authblk}

\usepackage{float}

\usepackage{booktabs,siunitx}
\usepackage{multirow}

\usepackage{fullpage}
\usepackage{bm}

\usepackage[super]{nth}

\usepackage{algpseudocode, algorithm, algorithmicx}

\usepackage{setspace}

\newcommand\amsclass[1]%
  {\hspace{9mm}
   \textbf{AMS subject classifications. }#1
}

\usepackage[T1]{fontenc}
\usepackage{romannum}
\usepackage{lineno}

\title{Stochastic approach for elliptic problems in perforated domains}
\date{}
\author{Jihun Han\thanks{jihun.han@dartmouth.edu}}
\author{Yoonsang Lee\thanks{yoonsang.lee@dartmouth.edu}}
\affil{Department of Mathematics, Dartmouth College}

\begin{document}
\pagenumbering{arabic}

\maketitle

\begin{abstract}
A wide range of applications in science and engineering involve a PDE model in a domain with perforations, such as perforated metals or air filters. Solving such perforated domain problems suffers from computational challenges related to resolving the scale imposed by the geometries of perforations. We propose a neural network-based mesh-free approach for perforated domain problems. The method is robust and efficient in capturing various configuration scales, including the averaged macroscopic behavior of the solution that involves a multiscale nature induced by small perforations. The new approach incorporates the derivative-free loss method that uses a stochastic representation or the Feynman-Kac formulation. In particular, we implement the Neumann boundary condition for the derivative-free loss method to handle the interface between the domain and perforations. A suite of stringent numerical tests is provided to support the proposed method's efficacy in handling various perforation scales. 
\end{abstract}



\section{Introduction}
Perforated materials play essential roles in science and engineering. Perforated metals in construction or mechanical engineering have significantly reduced weight while maintaining desired thermal and mechanical properties \cite{perfreduceweight}. Also, perforation can dampen sound waves to control noise for comfort and safety \cite{perfnoisereduction}. Electronics also use perforated materials to shield components from interference while achieving sufficient ventilation and preventing overheating \cite{perfwave}. Other uses of perforated materials include reducing light scattering \cite{perflightscattering}. In addition to the various uses of perforated materials, the perforation enables sustainable material design through the secondary use of materials, reducing harmful emissions in the melting process and prolonging the lifecycle.

Perforated materials or fluid systems involve processes in multiple spatiotemporal scales; the relatively small size of inclusions compared to the whole domain induces a wide range of scales to resolve. The dynamical scales of the sea ice with salt range from millimeters to tens of thousands of kilometers in space \cite{seaice}. In the battery pack case, the scale ratio between the large-scale battery pack and the cell is more than two orders of magnitude \cite{liionbattery}. Due to the wide range of spatiotemporal scales of perforated systems, it is computationally intractable to resolve all relevant scales. To make an accurate and robust prediction for such systems, it is essential to develop reliable and robust solvers for efficient representation of multiscale characteristics in the system, which are amenable to the current and the next generation computing resources.

Standard numerical methods are challenged in resolving all relevant scales due to their tremendously high computational costs. For heterogeneous media with specific structures, such as periodicity or ergodicity, the homogenization theory provides a closure model that describes macroscopic variations of the multiscale solution \cite{BLP}. In perforated domain problems, there are homogenized models under certain restrictions on the shapes and sizes of the inclusion, along with their relative distances (see \cite{PoissonRandomPerforation} for the Poisson equation, and \cite{StokesPerforatedSphere,  incompressible} for the Stokes equation and the derivation of Darcy's flow \cite{DarcyPerforated}). However, their applications are limited to special shapes (such as spheres) under the scale separation assumption.

There are research efforts in developing reduced-order numerical solvers for perforated domains, including the generalized multiscale finite element method (GMsFEM; \cite{GMsFEM}). GMsFEM involves two stages: offline and online. In the offline stage, a reduced-order space is constructed, which is then used to build multiscale basis functions resolving all fine structure details, such as variations related to small perforations, in the online stage. The key idea of such a two-stage process is to design appropriate snapshot spaces and determine an appropriate local spectral problem to select important modes in the snapshot space. This idea has been successfully applied to perforated domain problems, including various boundary conditions and weak PDE formulations (for example, mixed FEM) \cite{GMsFEM_eric,GMsFEM_robinBDC,MixedGMsFEM_eric,MsFEMAdvDifPerforated}. The computational efficiency of such a method is significantly obtained when the system is solved several times for various input parameters (ex, different forcing terms). In other words, GMsFEM still has a linear complexity in resolving all relevant scales to construct the multiscale basis functions.

The overarching goal of the current study is efficient sub-linear complexity numerical methods for perforated domains. The basis of the proposed approach is the derivative-free loss method (DFLM; \cite{DFLM}). DFLM uses a stochastic representation of a specific class of problems (see Eq.~\eqref{eq:model} in the next section), which involves an intrinsic averaging process to capture macroscopic behaviors. In other words, the solution is represented as an expectation of a martingale driven by a stochastic process, also known as the Feynman-Kac formulation. DFLM utilizes local stochastic walkers to explore neighborhoods for each collocation point, aligning with the homogenization rationale. 
For the solution function space, DLFM uses a neural network as a mesh-free representation of the solution. One characteristic of DFLM is that it differs from other neural network-based methods that use a point-wise residual loss function where the neighborhood variations are learned passively through a neural network. 
DFLM has recently been applied to elliptic homogenization in a simply connected domain without perforations \cite{DFLMHomo}, where the multiscale nature comes from the diffusion coefficient. In addition to the classical periodic structure case, the method successfully handled non-separable scale random field coefficients. The technique captures the homogenized solution without knowing the homogenized coefficient or calculating it through cell problems.

In formulating perforated domain problems, the perforation boundary conditions must be specified for the well-posedness. As the perforation does not allow any physical movement or fluctuations between the domain and perforation, the perforation boundary condition is often described as the homogeneous Neuman boundary condition. Under sufficient regularity assumptions on the boundary, we use the reflected stochastic process \cite{reflectivestochasticprocess} to implement the Neumann boundary condition, followed by the symmetrized Euler scheme for a numerical implementation \cite{symmetrizedEuler}.

We will first formulate the model problem of interest in Section \ref{sec:setup} and discuss some homogenization results for particular perforation configurations. The general description of the standard DFLM will be in Section \ref{sec:DFLM}, along with the Dirichlet boundary treatment.
Our approach for implementing the Neumann boundary condition Section is described in detail in \ref{sec:new}, where we also discuss time stepping to resolve the scales imposed by inclusion geometries. Section \ref{sec:experiment} will provide a suite of numerical results supporting the efficacy of DFLM in handling various perforation configurations. We will finish the paper with discussions of limitations and future directions of the proposed work in 
Section \ref{sec:discussion}.

%
%
%
\section{Problem setup}\label{sec:setup}
This section discusses the model problem of interest, an elliptic PDE in a perforated domain, and the homogenization result under a periodic structure, which provides an analytic model for the global scale behavior of the solution without perforations. Although its application of the homogenization theory is limited, we provide the theoretical result to validate our new method in Section \ref{sec:experiment}.

\begin{figure}[!t]
\centering
\includegraphics[width=0.5\textwidth]{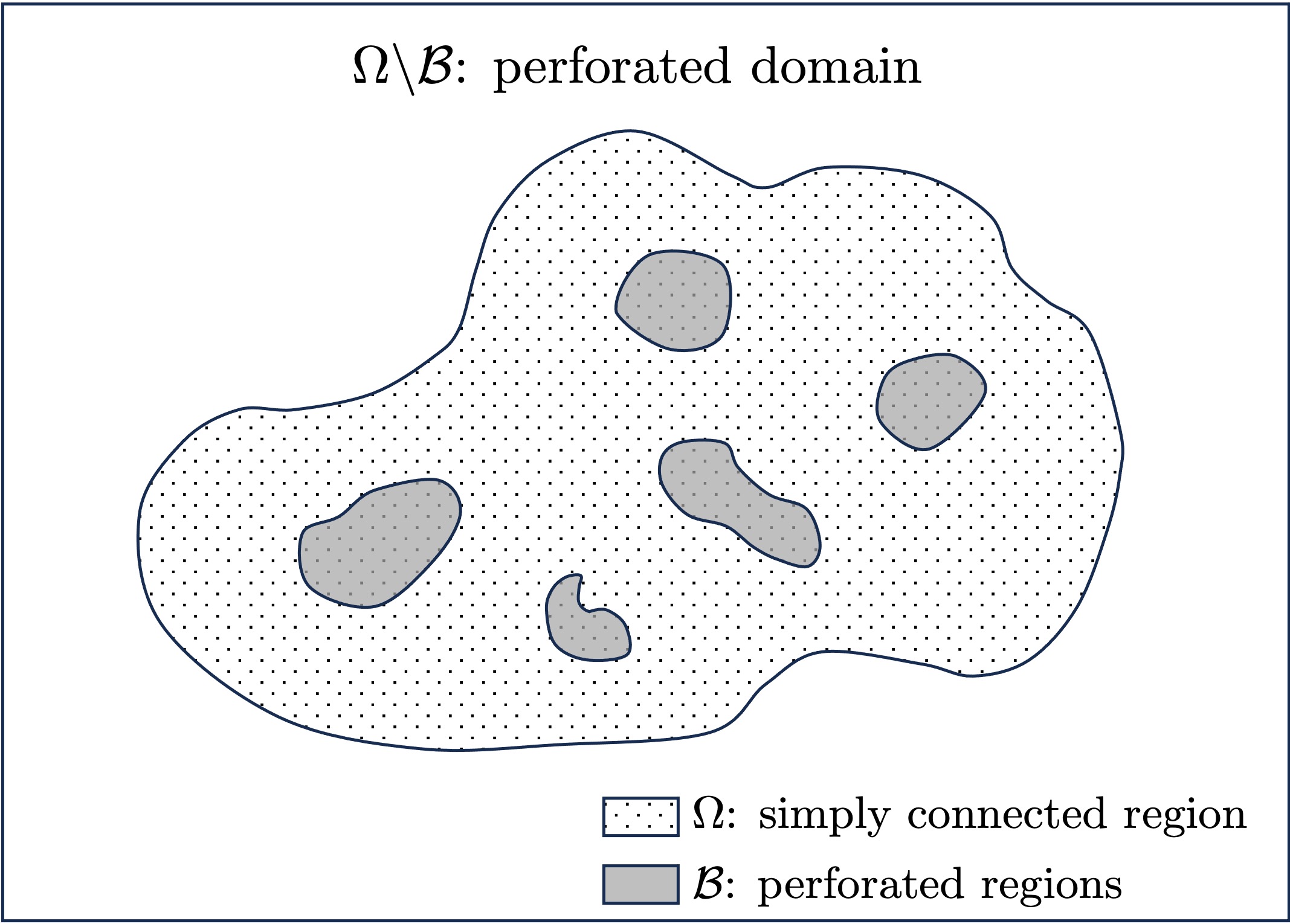}
\caption{A schematic of a perforated domain}\label{fig:perforation:schematic}
\end{figure}
For a collection of smooth perforations, $\mathcal{B}:=\cup_{i}\mathcal{B}_i$, and a simply connected domain $\Omega$ with smooth boundary (see a schematic of the perforated domain in Figure \ref{fig:perforation:schematic}), we consider an elliptic PDE in the perforated domain $\Omega\backslash\mathcal{B}$
\begin{equation}\label{eq:model}
\begin{split}
\frac{1}{2}\Delta u+\bm{V}(\bm{x},u)\cdot\nabla u&={G}(\bm{x},u)\quad\mbox{in}\quad \Omega\backslash\mathcal{B},\\
u&=g\quad\mbox{on}\quad\partial \Omega,\\
\frac{\partial u}{\partial \bm{n}}&=0\quad\mbox{on}\quad\partial \mathcal{B}.
\end{split}
\end{equation}
Here $\bm{V}=\bm{V}(\bm{x},u(\bm{x})) \in \mathbb{R}^k$ is the advection velocity field and $G=G(\bm{x}, u(\bm{x})) \in \mathbb{R}$ is the force term, both of which can depend on $u$. While the boundary condition on $\partial \Omega$ can take other types, such as Neumann or Robin, the boundary condition on the perforation is typically given by the Neumann condition, as it is challenging to know (or directly measure) the values on the perforation boundaries, $\partial \mathcal{B}$, to formulate the problem. In solving such a perforated domain problem numerically, one of the challenging issues comes from the discretization of the domain to account for different perforation configurations. 

When the scale of the perforation is relatively small compared to the global domain size, such a perforated domain problem involves a wide range of scales, and thus, resolving all relevant scales will be challenging. If there is a specific structure for the perforation, such as periodicity, a homogenization theory provides a relatively simple model to solve. In the current study, we consider the homogenization of the Poisson equation with the homogeneous boundary condition,
\begin{equation}\label{eq:Poissonmodel}
\begin{split}
\Delta u &= G\quad\mbox{in}\quad\Omega\backslash \mathcal{B},\\
u&=g\quad\mbox{ on } \quad\partial \Omega,\\
\frac{\partial u}{\partial \bm{n}}&=0\quad \mbox{ on } \quad \partial\mathcal{B}.
\end{split}
\end{equation}
Each perforation $\mathcal{B}_i$ is assumed to be rescaled from a perforation $\mathcal{B}_0$ in the hypercube $\square$ in $\mathbb{R}^d$
\begin{equation}
\square=\{x\in\mathbb{R}^d, |x_i|\leq 1/2\}
\end{equation}
with a ratio $\epsilon\ll 1$. When the scaling parameter $\epsilon$ approaches 0, the solution $u$ of Eq.~\eqref{eq:Poissonmodel} weakly converges to $u^0$ in $H^1$ where $u^0$ is the solution of the homogenized equation \cite{jikov, PoissonRandomPerforation}
\begin{equation}\label{eq:homogenized}
\begin{split}
\nabla A^0\nabla u^0 &= G \quad\mbox{in} \quad \Omega,\\
u^0&=g\quad\mbox{on}\quad\partial \Omega.\\
\end{split}
\end{equation}
Note that the homogenized equation does not involve the perforation, which is more straightforward to solve than the original perforated domain problem. Here, $A^0$ is the homogenized tensor, which can be a full matrix based on the geometries of the perforation, even for the Poisson problem. In particular, the homogenized tensor can be calculated by 
\begin{equation}\label{eq:homogenizedcoeff}
\langle \chi(\xi+\nabla N)\rangle=A^0\xi
\end{equation}
where $\chi$ is an arbitrary vector in $\mathbb{R}^d$ and $N\in H^1(\square\backslash \mathcal{B}_0)$ is a solution to the cell problem
\begin{equation}\label{eq:cell}
\begin{split}
\nabla(\chi(\xi+\nabla N))&=0 \quad\mbox{in}\quad \square\backslash \mathcal{B}_0,\\
\frac{\partial N}{\partial \bm{n}} \bigg|_{\partial \mathcal{B}_0}&=-\bm{n}\cdot\xi.
\end{split}
\end{equation}
For $\mathcal{R}^d$, by choosing a set of linearly independent $d$ $\chi$'s, we can find $A^0$ by inverting from \eqref{eq:homogenizedcoeff}.

If there is no specific structure in the perforation or scale separation (that is, no gap between the $\epsilon$ and $\Omega$ scales), no such homogenization model exists. To capture the macroscopic part of the multiscale solution, we incorporate a mesh-free approach based on neural networks and stochastic representation of the PDE Eq.~\eqref{eq:model}, described in detail in the following two sections. 

%
%
%
\section{Derivative-free loss method}\label{sec:DFLM}
We use a neural network as a mesh-free representation of multiscale solutions to capture the effective solution behavior of perforated domain problems with various perforation scales. In particular, we use the derivative-free loss method (DFLM) \cite{DFLM} that involves an intrinsic averaging representation through a stochastic representation of certain type PDEs. In this section, we describe the standard version of DFLM to approximate the model equation in a simply connected domain $\Omega$ without perforations
\begin{equation}\label{eq:modelOmega}
\begin{split}
\frac{1}{2}\Delta u+\bm{V}(\bm{x},u)\cdot\nabla u&={G}(\bm{x},u) \quad\mbox{in}\quad \Omega.\\
\end{split}
\end{equation}
Neumann boundary treatment, along with the Dirichlet case, will be discussed in the next section.

DFLM tackles the PDEs in the form of Eq.~\eqref{eq:modelOmega} through the underlying stochastic representation in the spirit of the Feynman-Kac formula that describes how the solution at a point interrelates to its neighborhood. DFLM guides a neural network directly to satisfy the interrelation across the entire domain, leading to the solution of the PDE. This intercorrelation characteristic makes DFLM different from other network-based methods that use pointwise residuals of the strong form PDE, such as PINN \cite{PINN}, in which a neural network implicitly learns the interrelation among different points.
Moreover, DFLM gradually and alternately updates a neural network and corresponding target values toward the PDE solution in the same manner as bootstrapping in the context of reinforcement learning, which differs from supervised learning methods that optimize neural network parameters within a fixed topology defined by loss functions.

For the stochastic representation of the solution of Eq.~\eqref{eq:modelOmega}, we consider the stochastic process $q({t};u, \bm{x}, \{\bm{X}_s\}_{0\leq s\leq t}) \in \mathbb{R}$ defined as 
\begin{equation}\label{eq:q_martingale_stochastic_process}
q({t};u, \bm{x},\{\bm{X}_s\}_{0\leq s\leq {t}}) := u(\bm{X}_{\Delta t}) - \int_0^{t} G(\bm{X}_s,u(\bm{X}_s))ds,\\
\end{equation}
where $\bm{X}_t \in \mathbb{R}^k$ is a stochastic process of the following SDE
\begin{equation}\label{eq:q_stochastic_walkers}
d\bm{X}_t = \bm{V}(\bm{X}_t, u(\bm{X}_t))dt + d\bm{B}_t,\quad \bm{B}_t:\mbox{ standard Brownian motion in }\mathbb{R}^k.
\end{equation}
Using It$\hat{\text{o}}$'s lemma (e.g., in \cite{karatzas2012brownian}), we find that the process $q({t};u, \bm{x}, \{\bm{X}_s\}_{0\leq s\leq t})$ satisfies the martingale property, which leads the representation of the solution $u(\bm{x})$ as follows;
\begin{align}\label{eq:q-martingale}
\begin{split}
u(\bm{x})&=q(0;u,\bm{x}, \bm{X}_0)=\mathbb{E}\left[q\left({\Delta t};u,\bm{x}, \{\bm{X}_s\}_{0\leq s\leq {\Delta t}} \right) | \bm{X}_0=\bm{x}\right] \\
&= \mathbb{E}\left[u(\bm{X}_{\Delta t}) - \int_{0}^{\Delta t}G(\bm{X}_s, u(\bm{X}_s)) ds \middle |\bm{X}_0=\bm{x}  \right],~\forall\bm{x}\in \Omega, \forall {\Delta t}>0.
\end{split}
\end{align}
Regarding to the definition of the stochastic process $q({t};u, \bm{x},\{\bm{X}_s\}_{0\leq s\leq {t}})$, the infinitesimal drift of the stochastic process $u(\bm{X}_t)$ is connected to the differential operator $\mathcal{N}[u]$ as 
\begin{equation}\label{eq:u_x_ito}
d(u(\bm{X}_t)) = \left(\mathcal{N}[u](\bm{X}_t) + G(\bm{X}_t, u(\bm{X}_t)\right) dt + \nabla u (\bm{X}_t) \cdot d\bm{B}_t.
\end{equation}
The martingale property, Eq.~\eqref{eq:q-martingale}, shows that the solution at a point $\bm{x}$, $u(\bm{X})$ can be represented through its neighborhood statistics observed by the stochastic process $\bm{X}_t$ starting at the point $\bm{x}$ during the time period $[0, {\Delta t}]$. That is, DFLM intrinsically involves averaging to represent the solution. 

DFLM constructs the loss function for training a neural network as
\begin{align}\label{eq:q_loss_continuous}
\mathcal{L}^{\Omega}(\bm{\theta})&=\mathbb{E}_{\bm{x}\sim \Omega}\left[\left|u(\bm{x};\bm{\theta})-\mathbb{E}\left[q\left(\Delta t;u(\cdot;\bm{\theta}),\bm{x}, \{\bm{X}_s\}_{0\leq s\leq \Delta t}\right) | \bm{X}_0=\bm{x}\right] \right|^2\right]\\
&=\mathbb{E}_{\bm{x}\sim \Omega} \left[\left|u(\bm{x};\bm{\theta}) - \mathbb{E}_{\{\bm{X}_s\}_{0\leq s\leq \Delta t}}\left[u(\bm{X}_t;\bm{\theta}) - \int_{0}^{\Delta t}G(\bm{X}_s, u(\bm{X}_s;\bm{\theta})) ds \middle |\bm{X}_0=\bm{x}  \right] \right|^2 \right]
\end{align} 
where the outer expection is over the sample collocation point $\bm{x}$ in the domain $\Omega$ and the inner expectation is over the stochastic path $\bm{X}_t$ starting at  $\bm{X}_0=\bm{x}$ during $[0,\Delta t]$. 
Analytically, the above formulation holds for an arbitrary macro time step  $\Delta t>0$. However, in the context of trainability, it is required to use a sufficiently long period to see the variation of the solution \cite{DFLManalysis}. If $\Delta t$ is too short, the distribution of $\bm{X}_t$ is close to the Dirac-Delta measure. Thus, the martingale representation Eq.~\eqref{eq:q-martingale} holds for an arbitrary initial guess, which hinders the training based on the $q$-martingale loss Eq.~\eqref{eq:q_loss_continuous}.

To calculate the stochastic path integral, $\bm{X}_t$ must be sampled at $\{m\delta t\}, m\in\mathbb{N}$ with a micro time step $\delta t$. In particular, with the presence of $u$ in the drift term $\bm{V}$, the computation of $\bm{X}_t$ by solving the stochastic differential equation can be costly. As an alternative to avoid the calculation of the solution to Eq.~\eqref{eq:q_stochastic_walkers}, another martingale process $\tilde{q}({\Delta t};u,\bm{x},\{\bm{B}_s\}_{0\leq s\leq {\Delta t}})$ based on the standard Brownian motion $\bm{B}_t$ is proposed as 
\begin{align}\label{eq:q_tilde_martingale_process}
\tilde{q}({\Delta t};u, \bm{x}&,\{\bm{B}_s\}_{0\leq s\leq {\Delta t}}) := 
 \left(u(\bm{B}_{\Delta t}) - \int_0^{\Delta t} G(\bm{B}_s,u(\bm{B}_s))ds\right) \mathcal{D}(\bm{V},u,{\Delta t}),\\
&
\text{where}{\quad}\mathcal{D}(\bm{V},u,{\Delta t})=\exp \biggl(\int^{\Delta t}_0\bm{V}(\bm{B}_s, u(\bm{B}_s))\cdot d\bm{B}_s-\frac{1}{2}\int^{\Delta t}_{0}|\bm{V}(\bm{B}_s, u(\bm{B}_s))|^2 ds \biggr). \nonumber 
\end{align}
Here, $\tilde{q}$-process replaces $\bm{X}_t$ to $\bm{B}_t$ in the $q$-process with an additional exponential factor $\mathcal{D}(\bm{V}, u,t)$ compensating the removal of the drift effect in $\bm{X}_t$ \cite{karatzas2012brownian, oksendal}. The alternative $\tilde{q}$-martingale allows the standard Brownian walkers to explore the domain regardless of the form of the given PDE, which can be drawn from the standard Gaussian distribution without solving SDEs, thus saving a significant amount of computational cost. Using the computationally efficient $\tilde{q}$-martingale, the loss function is 
\begin{equation}
\mathcal{L}^{\Omega}(\bm{\theta})=\mathbb{E}_{\bm{x}\sim \Omega}\left[\left|u(\bm{x};\bm{\theta})-\mathbb{E}\left[\tilde{q}({\Delta t};u(\cdot;\bm{\theta}),\bm{x}, \{\bm{B}_s\}_{0\leq s\leq {\Delta t}}) | \bm{B}_0=\bm{x}\right] \right|^2\right].
\end{equation}


The loss function ${\mathcal{L}}^{\Omega}(\bm{\theta})$ is optimized by a stochastic gradient descent method, and, in particular, the bootstrapping approach is used as the target of the neural network (i.e., the expectation component of $q$- or $\tilde{q}$-process) is pre-evaluated using the current state of neural network parameters $\bm{\theta}$. The $n$-th iteration step for updating the parameters $\bm{\theta}_n$ is 
\begin{equation}\label{eq:sgd_update}
\bm{\theta}_n = \bm{\theta}_{n-1} - \alpha \nabla{\mathcal{L}}_n^{\Omega}(\bm{\theta}_{n-1}).
\end{equation}
Here, ${\mathcal{L}}_n^{\Omega}(\bm{\theta})$ is the empirical interior loss function using $N_r$ sample collocation points $\{\bm{x}_i\}_{i=1}^{N_r}$ in the interior of the domain $\Omega$, and $N_s$ stochastic walkers $\left\{\bm{X}_s^{(i,j)};s \in [0,\Delta t],  \bm{X}_0=\bm{x}_i\right \}_{j=1}^{N_s}$ at each sample collocation point $\bm{x}_i$, 
\begin{align}\label{eq:q_loss_discrete}
{\mathcal{L}}^{\Omega}_{n}(\bm{\theta}) :&= {\mathbb{E}}_{\bm{x}\sim \Omega}\left[\left|u(\bm{x};\bm{\theta})-{\mathbb{E}}\left[q\left(\Delta t;u(\cdot;\bm{\theta}_{n-1}),\bm{x}, \{\bm{X}_s\}_{0\leq s\leq \Delta t}\right) | \bm{X}_0=\bm{x}\right] \right|^2\right] \\
&= \frac{1}{N_r} \sum \limits_{i=1}^{N_r}\left| u(\bm{x}_i;\bm{\theta})-\frac{1}{N_s}\sum\limits_{j=1}^{N_s}\left\{ u\left(\bm{X}_{\Delta t}^{(i,j)};\bm{\theta}_{n-1} \right)-\int_0^{\Delta t}G\left(\bm{X}_s^{(i,j)}, u\left(\bm{X}_{\Delta t}^{(i,j)};\bm{\theta}_{n-1} \right)\right) ds \right\} \right |^2.
\end{align}
The random interior collocation points $\{\bm{x}_i\}_{i=1}^{N_r}$ can follow a distribution whose support covers the domain $\Omega$. The learning rate $\alpha$ could be tuned at each step, and the gradient descent step can be optimized by considering the previous steps, such as Adam optimization \cite{adam}.

\section{DFLM for perforated domains}\label{sec:new}
The standard DFLM \cite{DFLM} (and its extension for fluid problems \cite{DFLMfluid}) has been implemented for the Dirichlet boundary condition. As the typical conditions on the perforation boundaries are homogeneous Neumann, we propose how to implement homogeneous Neumann boundary conditions in this section.  DFLM implicitly implements the boundary condition by correctly handling the stochastic process near the boundaries using killed (Dirichlet) or reflected (Neumann) diffusion, which guarantees a martingale representation. In other words, DFLM avoids an additional loss term related to the boundary conditions and thus bypasses the balancing between different loss components, which can affect the training performance. 

\subsection{Dirichlet Boundary}
To impose the Dirichlet boundary condition for the stochastic representation of the PDE solution, we introduce a stopping time with respect to the process $\bm{X}_t$ defined as $\tau:= \inf \left\{s: {\bm{X}}_s \notin \Omega \right \}$. It is the first hitting time of the process to the boundary $\partial \Omega$. Then the stopped process ${\bm{X}}_{\min \{{\Delta t}, \tau \}}$ indicates that the process is absorbed by the boundary $\partial \Omega$ upon reaching it. From the optimal stopping theorem \cite{karatzas2012brownian}, the stopped $q$-process, $q(\min\{\Delta t, \tau\};u,\bm{x},\{\bm{X}_s\}_{0\leq s \leq \min\{\Delta t, \tau\}})$, is also martingale, and we have the representation of the PDE solution
\begin{equation}\label{eq:q-dirichelt-neumann-martingale}
u(\bm{x})=\mathbb{E}\left[u({\bm{X}}_{\min \{\Delta t,\tau \}}) - \int_{0}^{\min \{\Delta t,\tau \}}G({\bm{X}}_s, u({\bm{X}}_s)) ds \middle |{\bm{X}}_0=\bm{x}  \right],
\end{equation}
which generalizes \eqref{eq:q-martingale}. We note that Eq.~\eqref{eq:q-dirichelt-neumann-martingale} also holds for $\bm{x}\in \partial \Omega$ since the stopped $q$-process on the right-hand side is stationary with its initial value $u(\bm{x})$ being identical to the prescribed value $g(\bm{x})$. Moreover, when stochastic paths $\bm{X}_t$ with $\bm{X}_0 = \bm{x} \in \Omega$ hit the boundary $\partial \Omega$ within the time period $[0,\Delta t]$ (i.e., $\min\{\Delta t, \tau\}=\tau$), we effectively transfer the exact boundary value information to the interior point $\bm{x}$ by replacing $u({\bm{X}}_{\min \{\Delta t,\tau \}})$ with $g({\bm{X}}_{\min \{\Delta t,\tau \}})$ in the expectation computation on right-hand side. In sum, \eqref{eq:q-dirichelt-neumann-martingale} accomodates both i) the PDE solution inside the domain and ii) Dirichlet condition on $\partial \Omega$. Similarly, we can derive a comparable representation utilizing a stopped $\tilde{q}$-process driven by the standard Brownian motion $\bm{B}_t$ incorporating the stopping time $\tau^{\prime}:=\inf \{s:\bm{B}_s \notin \Omega\}$;
\begin{equation}\label{eq:stoppted_q_tilde_representation}
u(\bm{x}) = \mathbb{E}\left[\tilde{q}({\min \{\Delta t, \tau^{\prime}\} };u(\cdot;\bm{\theta}),\bm{x}, \{\bm{B}_s\}_{0\leq s\leq {\min \{\Delta t, \tau^{\prime}\}}}) | \bm{B}_0=\bm{x}\right].
\end{equation}

\begin{figure}[t!]
\captionsetup[subfigure]{justification=centering}
\begin{subfigure}[b]{0.49\textwidth}
\centering
\includegraphics[width=0.85\textwidth]{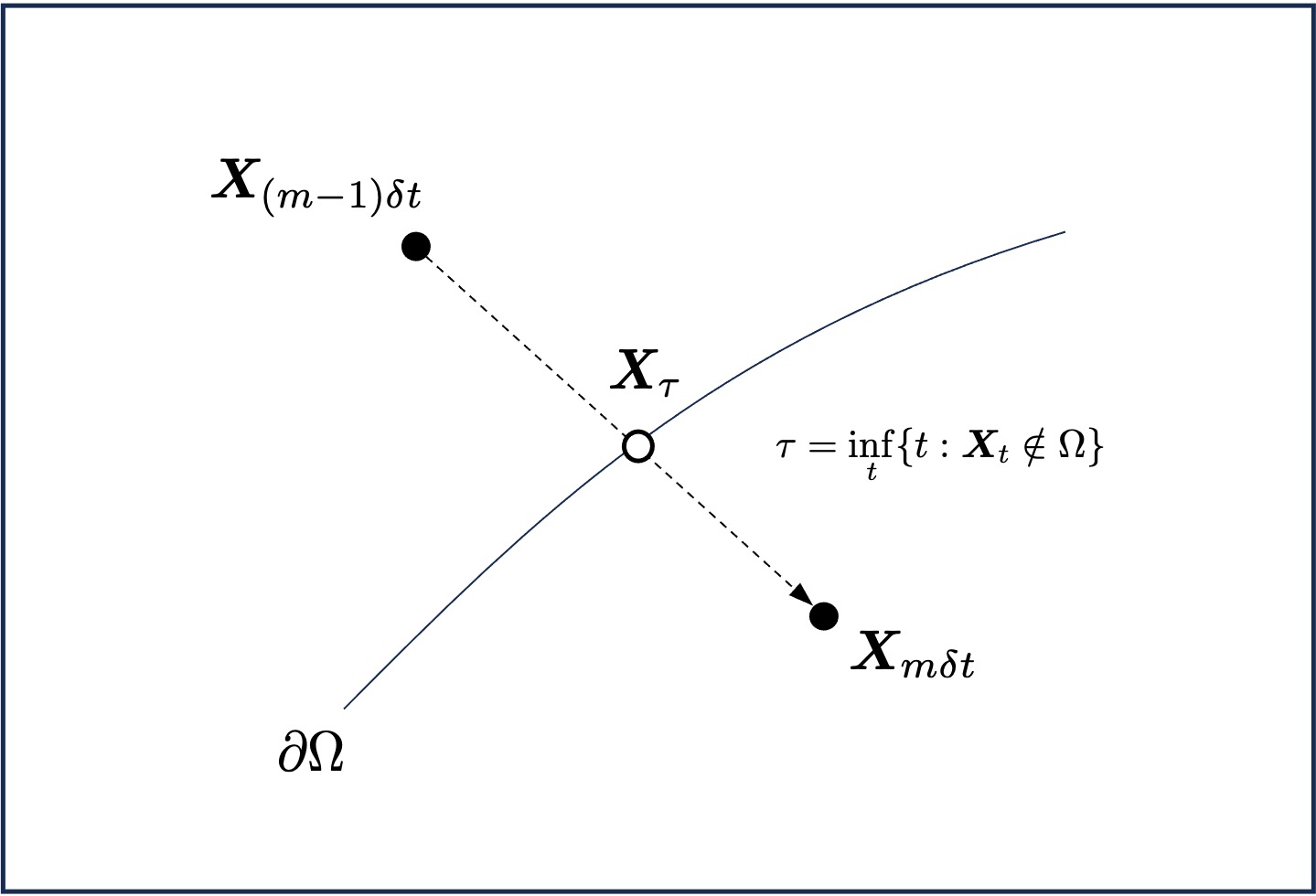}
\caption{Dirichlet boundary condition \\ (killed diffusion)}
\end{subfigure}
\begin{subfigure}[b]{0.49\textwidth}
\centering
\includegraphics[width=0.85\textwidth]{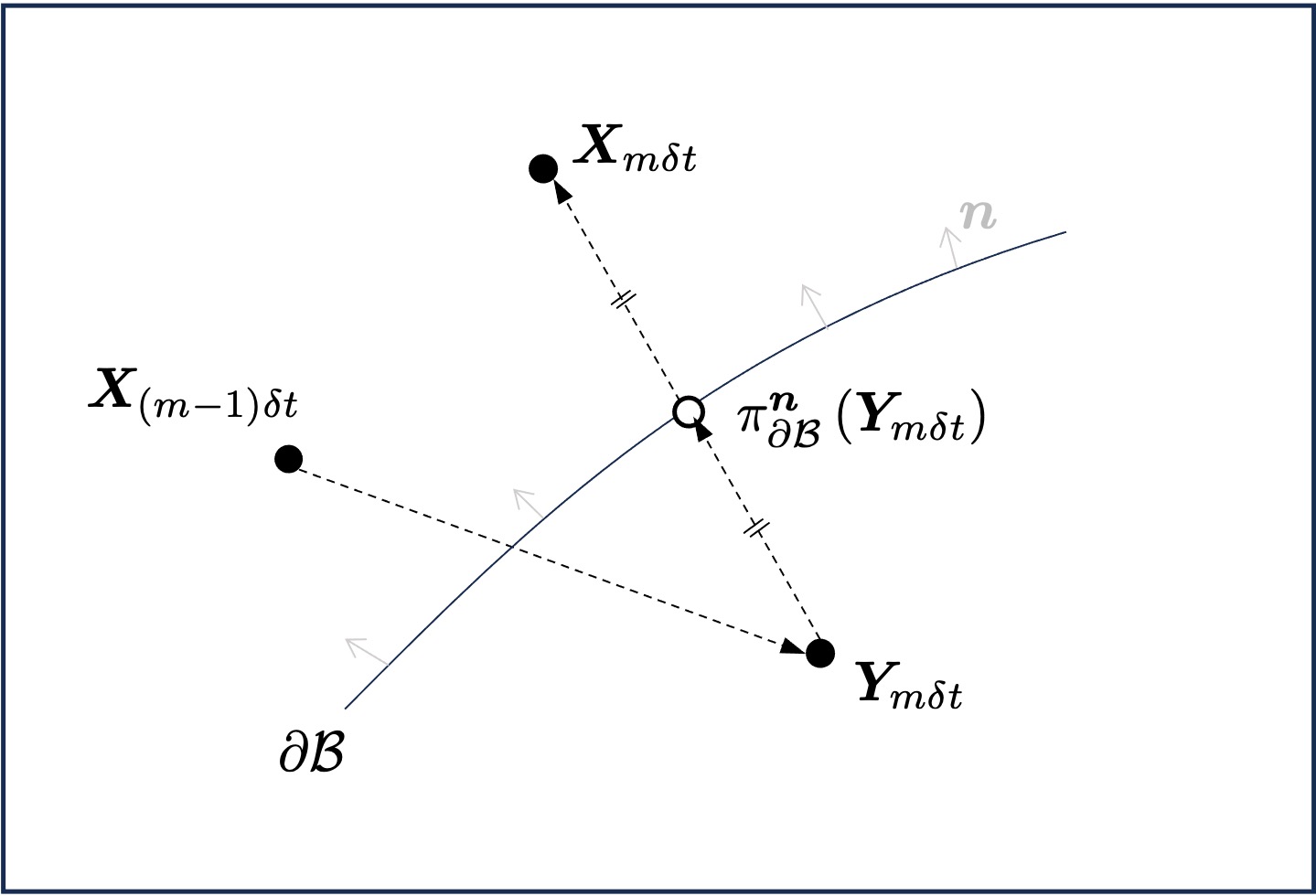}
\caption{Neumann boundary condition \\  (reflected diffusion)}
\end{subfigure}
\caption{Schematic diagram of numerical treatment for boundary conditions.}
\label{fig:BC_diagram}
\end{figure}

The crucial part is to estimate the accurate exit position. (In the subsequent discussion, we can also read $\bm{X}_t$ as $\bm{B}_t$ and $\tau$ as $\tau^{\prime}$ in case $\tilde{q}$-representation Eq.~\eqref{eq:stoppted_q_tilde_representation}.) To resolve the dynamics of the stochastic process, we use a sufficiently small micro time step $\delta t$ and sample $\bm{X}_{t}$ at $\{j\delta t\}_{j=0}^M$ (where $M\delta t=\Delta t$).
When the process hits the boundary in $[(m-1)\delta t, m\delta t]$, we estimate the exit position $\bm{X}_{\tau}$ using the intersection of the line segment connecting ${\bm{X}}_{(m-1)\delta t}$ and ${\bm{X}}_{m\delta t}$ and the boundary, which implements the killed diffusion (see Figure \ref{fig:BC_diagram} (a) for a schematic of the killed diffusion). The exit time $\tau$ is linearly approximated using the length of walker's journey until termination given by  $\tau = (m-1+\alpha) \delta t$ where $\alpha=\frac{|\bm{X}_{\tau}-\bm{X}_{(m-1)\delta t}|}{|\bm{X}_{m\delta t}-\bm{X}_{(m-1)\delta t}|}$. 
Also, for scalar and vector values functions $h$ and $\bm{h$}, we compute the stochastic path integral using the linear approximation of the exit time;
\begin{equation}
\int_{(m-1)\delta t}^{m\delta t} h(\bm{X}_s)ds =\int_{(m-1)\delta t}^{\tau} h(\bm{X}_s)ds  ~~\simeq~~ h\left(\bm{X}_{(m-1)\delta t}\right)(\alpha \delta t),
\end{equation}
and, for the computation of additional exponential factor in Eq.~\eqref{eq:stoppted_q_tilde_representation} (see Eq.~\eqref{eq:q_tilde_martingale_process})
\begin{equation}
\int_{(m-1)\delta t}^{m\delta t} \bm{h}(\bm{B}_s)\cdot d\bm{B}_s = \int_{(m-1)\delta t}^{\tau} \bm{h}(\bm{B}_s)\cdot d\bm{B}_s ~~\simeq~~ \bm{h}(\bm{B}_{(m-1)\delta t})\cdot\left(\bm{B}_{\tau}-\bm{B}_{(m-1)\delta t}\right). 
\end{equation}

\subsection{Neumann Boundary}
With the physical interpretation of the no-flux for the homogeneous Neumann boundary condition, we consider the reflected stochastic process \cite{reflectivestochasticprocess} associated with the boundaries:
\begin{equation}\label{eq:q_stochastic_walkers_reflect}
d{\bm{X}}_t = \bm{V}({\bm{X}}_t, u({\bm{X}}_t)) dt + d\bm{B}_t + \bm{n}({\bm{X}}_t) dL_t.
\end{equation}
Here $\bm{n}$ is the inward normal vector on the boundary $\partial \mathcal{B}$, and $L_t$ is a local time process with $L_0=0$. $L_t$ is continuous and increases only when ${\bm{X}}_t$ is on the boundary $\partial \mathcal{B}$, defined as:
\begin{equation}
L_t = \int_0^{t} 1_{\{{\mathbf{X}}_s \in \partial \mathcal{B}\}} dL_s.
\end{equation}
We note that the stochastic process in Eq.~\eqref{eq:q_stochastic_walkers_reflect} stems from the modification of the stochastic process defined in Eq.~\eqref{eq:q_stochastic_walkers}, reflecting on the perforation boundaries. Such reflecting stochastic differential equation of $({\bm{X}}_t, L_t)$ as Eq.~\eqref{eq:q_stochastic_walkers_reflect} is commonly referred to as \textit{Skorokhod problem} or the \textit{problem of reflection} \cite{BCRegularity1,BCRegularity2}. The existence of the unique solution of Eq.~\eqref{eq:q_stochastic_walkers_reflect} is established under favorable regularity conditions on the domain \cite{BCRegularity3, BCRegularity4, BCRegularity5}. We assume that perforations have a sufficiently smooth boundary, such as class $C^5$, to guarantee the existence and uniqueness.

By It$\hat{\text{o}}$'s lemma,
\begin{align}
d\left(u({\bm{X}}_t)-\int_0^{\Delta t}G({\bm{X}}_s)ds\right) &= \mathcal{N}[u]({\bm{X}}_t) dt + \nabla u ({\bm{X}}_t) \cdot d\bm{B}_t + \nabla_{\bm{x}}u({\bm{X}}_t)\cdot \bm{n}({\bm{X}}_t) dL_t \\
& = \nabla u ({\bm{X}}_t) \cdot d\bm{B}_t
\end{align}
where the term $\nabla_{\bm{x}}u(\bm{X}_t)\cdot \bm{n}(\bm{X}_t) dL_t$ vanishes by the definition of the local  time process $L_t$ and the homogeneous Neumann boundary condition implies $\nabla_{\bm{x}}u(\bm{X}_t)\cdot \bm{n}(\bm{X}_t)= \frac{\partial u}{\partial n}(\bm{X}_t) =0$ on $\partial \mathcal{B}$.
\begin{figure}[h!]
\centering
\includegraphics[width=0.4\textwidth]{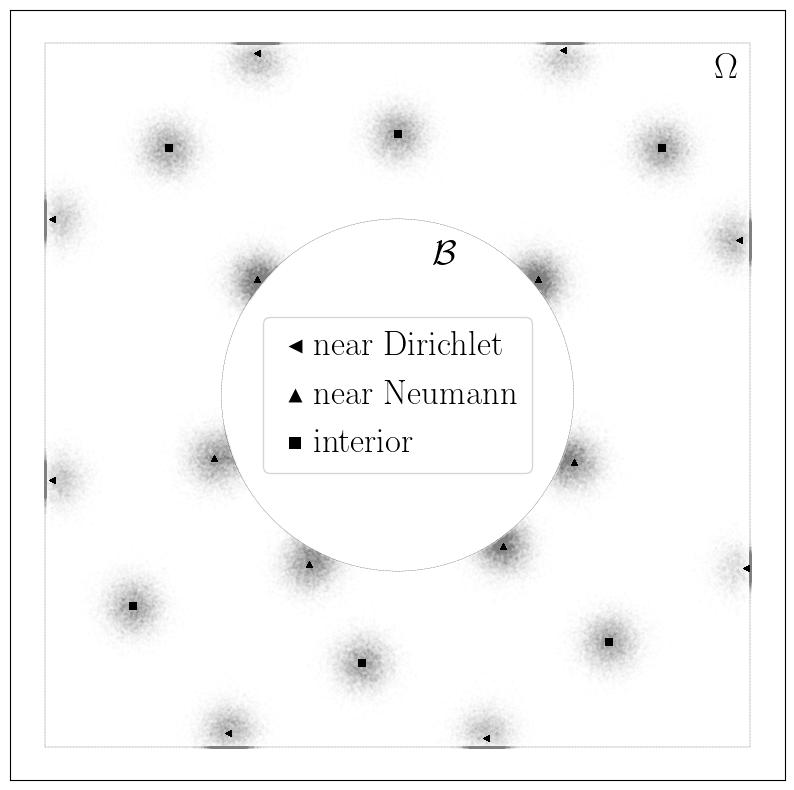}
\caption{Brownian sample distribution of colocation points based on the proximity to the boundary with Neumann and Dirichlet conditions on $\partial\mathcal{B}$ and $\partial \Omega$, respectively. Brownian samples are killed near $\partial \Omega$ and reflected near $\partial \mathcal{B}$, contrasting with the symmetric distribution of standard Brownian samples.}
\label{fig:BC_simulation}
\end{figure}
Therefore, the process, $u({\bm{X}}_t)-\int_0^{\Delta t}G({\bm{X}}_s)ds$, which is the same as $q$-process in Eq.~\eqref{eq:q-martingale} with  the modification of reflection Eq.~\eqref{eq:q_stochastic_walkers_reflect}, is also martingale. This leads to the representation of the PDE solution as Eq.~\eqref{eq:q-martingale}. Similarly, we can obtaion the represenation of PDE solution with reflected Brownian motion as $\tilde{q}$-martingale representation.

The numerical simulation of reflecting stochastic processes or the Skorohod problems uses either projection methods \cite{symmetrizedEuler, projectionNeumannSim1, projectionNeumannSim2, projectionNeumannSim3} or penalty methods \cite{penaltyNeumannSim1, penaltyNeumannSim2, penaltyNeumannSim3}. In this work, we use the symmetrized Euler scheme \cite{symmetrizedEuler} as a type of projection method to simulate the reflecting process Eq.~\eqref{eq:q_stochastic_walkers_reflect} associated with the Neumann boundary $\partial \mathcal{B}$. We evolve the stochastic process from $\bm{X}_{(m-1)\delta t}$, yielding $\bm{Y}_{m\delta t}$. If $\bm{Y}_{m\delta t}$ is in the perforation, reflect $\bm{Y}_{m\delta t}$ with respect to the perforation boundary. Otherwise, set $\bm{X}_{m\delta t}$ to be $\bm{Y}_{m\delta t}$. The scheme is summarized as follows:
\begin{align}
\text{(i)}&~~ \bm{Y}_{m \delta t} = {\bm{X}}_{(m-1) \delta t}+\bm{V}({\bm{X}}_{(m-1) \delta t},u({\bm{X}}_{(m-1) \delta t}))\delta t + \sqrt{\delta t} \bm{Z}, {\quad} \bm{Z} \sim \mathcal{N}(0, \mathbb{I}_k),\\
\text{(ii)}&~~ {\bm{X}}_{m \delta t} = 
\begin{cases}
\bm{Y}_{m\delta t} & ~~\textrm{if}~~~\bm{Y}_{m\delta t} \in \overline{\Omega\backslash \mathcal{B}} \\
\bm{Y}_{m\delta t} + 2\left(\pi_{\partial \mathcal{B}}^{\bm{n}}(\bm{Y}_{m \delta t})-\bm{Y}_{m \delta t}\right) & ~~ \textrm{if}~~~\bm{Y}_{m\delta t} \in \overline{\mathcal{B}}
\end{cases}
\end{align}
where $\pi_{\partial \mathcal{B}}^{\bm{n}}(\bm{Y}_{m \delta t})$ is the projection of $\bm{Y}_{m \delta t}$ onto $\partial \mathcal{B}$ parallel to $\bm{n}$. We describe the diagram for the scheme in Fig.~\ref{fig:BC_diagram}-(b). We note that our assumption on $\mathcal{B}$ with $C^5$-regularity guarantees the uniqueness of the projection \cite{symmetrizedEuler}. Also, in case of the reflected Brownian motion, you can modify the scheme by replacing $\bm{X}_t$ to $\bm{B}_t$ and setting $\bm{V}=0$.

When ${\bm{X}}_t$ reflects on the boundary $\partial \mathcal{B}$ during the time interval $[(m-1)\delta t, m \delta t]$, we estimate the reflection time $(m-1+\beta)\delta t$ through linear approximation
\[\beta = \frac{\left|\pi_{\partial D}^{\bm{n}}(\bm{Y}_{m \delta t})-\bm{Y}_{m \delta t} \right|}{\left| {\bm{X}}_{(m-1)\delta t} -\pi_{\partial D}^{\bm{n}}(\bm{Y}_{m \delta t}) \right|+\left|\pi_{\partial D}^{\bm{n}}(\bm{Y}_{m \delta t})-\bm{Y}_{m \delta t} \right|}.\]
Using the estimated $\beta$, the stochastic path integral for scalar and vector-valued functions $h$ and $\bm{h}$ are approximated through two path segments, one into the boundary and another out of the boundary
\begin{equation}
\int_{(m-1)\delta t}^{m\delta t} h({\bm{X}}_s)ds ~~\simeq~~ h\left({\bm{X}}_{(m-1)\delta t}\right)((1-\beta) \delta t) + h\left(\pi_{\partial \mathcal{B}}^{\bm{n}}(\bm{Y}_{m \delta t})\right)(\beta\delta t) 
\end{equation}
and, for the representation with reflecting Brownian motion, 
\begin{equation}
\int_{(m-1)\delta t}^{m\delta t} \bm{h}({\bm{B}}_s)\cdot d{\bm{B}}_s ~~\simeq~~ \bm{h}({\bm{B}}_{(m-1)\delta t})\cdot\left(\pi_{\partial \mathcal{B}}^{\bm{n}}(\bm{Y}_{m \delta t})-{\bm{B}}_{(m-1)\delta t}\right) + h(\pi_{\partial \mathcal{B}}^{\bm{n}}(\bm{Y}_{m \delta t}))\cdot(\pi_{\partial \mathcal{B}}^{\bm{n}}(\bm{Y}_{m \delta t})-\bm{Y}_{m \delta t}).
\end{equation}

In sum, we address both Dirichlet and Neumann boundary conditions by adjusting the characteristics of the stochastic processes $\bm{X}_t$ (or $\bm{B}_t$) in the vincity of the boundary either absorption for Dirichlet or reflection for Neumann condition. Fig.~\ref{fig:BC_simulation} presents the simulation example of stochastic samples at various colocation points in relation to their proximity to the boundary. As depicted, the sample distribution varies based on the specific boundary conditions applied.

\subsection{Time stepping}
The micro time step $\delta t$ needs to be sufficiently small to calculate the stochastic path integrals accurately. Another constraint on the micro time step is related to the perforation size. If $\delta t$ is too large while a perforation is small, there is a high chance of missing the perforation, and thus, the stochastic walker does not interact with the perforation boundary. We impose the following constraint on the micro time step $\delta t$
\begin{equation}
\mathbb{E}[|\bm{B}_{m\delta t}-\bm{B}_{(m-1)\delta t}|] \ll \mbox{perforation size}
\end{equation}
to ensure that the line segment from $\bm{B}_{(k-1)\delta t}$ to $\bm{B}_{k\delta t}$ (approximating the Brownian trajectory within the micro-time step) intersects with the boundary of perforation at most once. Also, this constraint guarantees that the reflected process does not intersect with the other boundary. At most, one reflection happens during a micro time interval, almost surely. To be specific, from the properties of the Brownian motion, we have
\begin{equation}
\mathbb{E}[|\bm{B}_{m\delta t}-\bm{B}_{(m-1)\delta t}|]=\mathbb{E}[|\bm{B}_{\delta t}|]=\kappa(d)\sqrt{dt}
\end{equation}
where $\kappa(d)$ is a constant depending on the dimension of the domain. For example, in 2D, $\kappa(2)=\sqrt{\frac{\pi}{2}}$.
%
%
%
\section{Numerical Experiments}\label{sec:experiment}

In this section, we validate the effectiveness of DFLM in handling various perforations in the unit square $\Omega = [-0.5,0.5]^2 
\subset \mathbb{R}^2$. As the challenge comes from perforations, we focus on solving the Poisson problem in various perforation configurations with the following boundary conditions
\begin{align}
-\frac{1}{2}\Delta u &= 1 \quad \mbox{in} \quad  \Omega \backslash \mathcal{B}, \\
u &=1 \quad \mbox{on} \quad  \partial \Omega, \\
\frac{\partial u}{\partial \bm{n}}& = 0 \quad \mbox{on} \quad \partial \mathcal{B}.
\end{align}
As $G=0$ and $\bm{V}=0$, the corresponding martingale representation of the solution is given by 
\begin{equation}
u(\bm{x}) = \mathbb{E} \left[u(\bm{X}_{\Delta t}) + \Delta t  \middle | \bm{X}_0=\bm{x} \right], ~~\textrm{where}~~d\bm{X}_t = d \bm{B}_t + \bm{n}(\bm{X}_t)dL_t.
\end{equation}
In the current study, we consider circular perforations where the regularity of the boundaries ensures the existence of reflective Brownian motion. We denote each circular hole as $B(\bm{a};r):=\{\bm{x} : |\bm{x}-\bm{a}|\leq r\}$. 

In designing the network, we incorporate the Fourier-feature embedding neural network \cite{tancik2020fourier}. Note that the Fourier network is essentially a standard multilayer perceptron (MLP) with the inclusion of embedding input into Fourier space as follows:
\begin{equation}
\bm{x} \mapsto \left[\begin{matrix}\cos (\bm{A}_{\sigma}\bm{x}) \\ \sin (\bm{A}_{\sigma}\bm{x}) \end{matrix} \right] \in \mathbb{R}^{2m}, \quad \bm{A}_{\sigma} \in \mathbb{R}^{n\times m}, \quad (\bm{A}_{\sigma})_{ij} \sim \mathcal{N}(0,\sigma^2).
\end{equation}
This input transformation is known to attenuate the spectral bias during neural network training \cite{spectralBias}. 
We train neural networks in all experiments using stochastic gradient descent (SGD) with the Adam optimizer and learning parameters $\beta_1=0.99$ and $\beta_2=0.99$. 
All trainable parameters are initialized from the Glorot normal distribution, taking into account the dimensions of the input and output of each layer. We employ an exponential decay learning rate with an initial rate $\alpha_0$ and decay rate $\gamma$ per 1000 training iterations. The specific values of $\alpha_0$ and $\gamma$ will be specified for each test problem. Instead of updating $\bm{\theta}_{n}$ every step, we run three gradient descent steps using $\bm{\theta}_{n-1}$ in order to bring $u(\bm{x};\bm{\theta}_n)$ closer to the target $\mathbb{E}\left[\tilde{q}(t;u(\cdot;\bm{\theta}_{n-1}),\bm{x}, \{\bm{B}_s\}_{0\leq s\leq \Delta t}) | \bm{B}_0=\bm{x}\right]$:
\begin{equation}
\bm{\theta}_{n-1}^{(m+1)} = \bm{\theta}_{n-1}^{(m)} - \alpha \nabla {\mathcal{L}}_n\left(\bm{\theta}_{n-1}^{(m)}\right), \quad m=0,1,2, \quad \bm{\theta}_{n-1}^{(0)}=\bm{\theta}_{n-1}, \quad \bm{\theta}_{n-1}^{(3)} = \bm{\theta}_n.
\end{equation}
To evaluate the empirical loss Eq.~\eqref{eq:q_loss_discrete}, we sample $N_r$ training samples from a uniform distribution at each iteration using rejection sampling in a perforated domain. 

We obtain reference solutions $u$ using the finite element method (FEM) with sufficiently fine mesh sizes after ensuring convergence. The accuracy of a neural network solution $\tilde{u}$ is measured by the relative $\mathcal{L}^2$-error, defined as $\frac{\|u-\tilde{u}\|_{\Omega\backslash\mathcal{B}}}{\|u\|_{\Omega\backslash\mathcal{B}}}$, where the $\mathcal{L}^2$-norm is computed on equidistant $1001\times 1001$ grid points.

%
%
\subsection{Single perforation}
\begin{figure}[t!]
\centering
\includegraphics[width=1.0\textwidth]{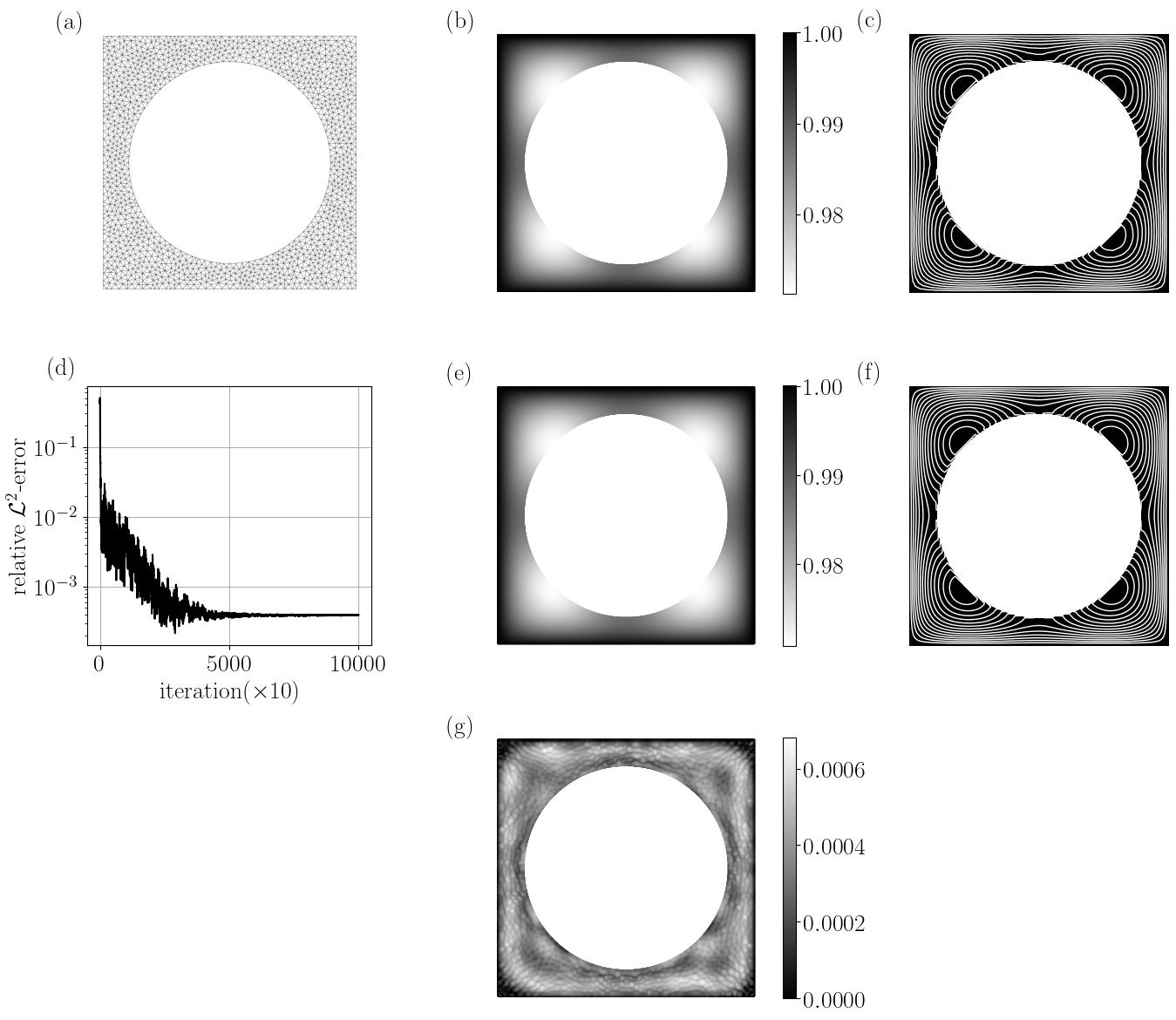}
\caption{Result of single perforation. (a): Finite element mesh, (b): FEM solution, (c): contour map of the FEM solution, (d): learning trajectory of DFLM, (e): DFLM approximation, (f): contour map of the DFLM solution, (g): pointwise difference.}
\label{fig:example1}
\end{figure}

The first example is a single perforation $\mathcal{B}=B((0,0);0.4)$, which is a circle of radius 0.4 centered at (0.0) (see Fig.~\ref{fig:example1} for the configuration). This single perforation problem is to verify the implementation of the Neuman boundary condition in DFLM through the reflected stochastic processes. 

The dimension of the Fourier embedding layer is 200 (i.e., $\bm{A}_{\sigma} \in \mathbb{R}^{100 \times 2}$) with the variance $\sigma^2$ set to 9. Three hidden layers of dimension 200 follow the first layer with $\tanh$ activation. $ N_r=3000$ center collocation points evaluate the empirical training loss with $N_s=400$ stochastic walkers at each collocation point. The learning rate parameters $(\alpha_0,\gamma)$ are $(10^{-3}, 0.9)$. In this example, the perforation size is relatively large (comparable to the size of $\Omega$) and thus, the micro step $\delta t$ is chosen mainly considering the accuracy of the stochastic integral, which is $5 \times 10^{-6}$ while the macro time step $\Delta t$ is $128\times \delta t$.

Fig.~\ref{fig:example1} shows the DFLM result (second row) along with the reference FEM solution (first row). FEM requires a relatively small discretization to resolve the circular shape of the perforation to achieve convergence (Fig.~\ref{fig:example1}-(a)). The DFLM solution converges after 60,000 iterations, which is shown in Fig.~\ref{fig:example1}-(e), along with its contour lines (Fig.~\ref{fig:example1}-(f)). The DFLM contour lines are perpendicular to the boundary, which verifies that the normal derivative of the solution vanishes on the interface, showing correct implementation of the homogeneous Neumann boundary condition. In comparison with the reference solution, the DFLM achieves pointwise errors of order $5\times 10^{-4}$ (see Fig.~\ref{fig:example1}-(g) for the pointwise errors). The relative $\mathcal{L}^2$-error is $3.9023 \times 10^{-4}$.

%
%
\subsection{Perforation of various scales}
\begin{figure}[!t]
\centering
\includegraphics[width=1.0\textwidth]{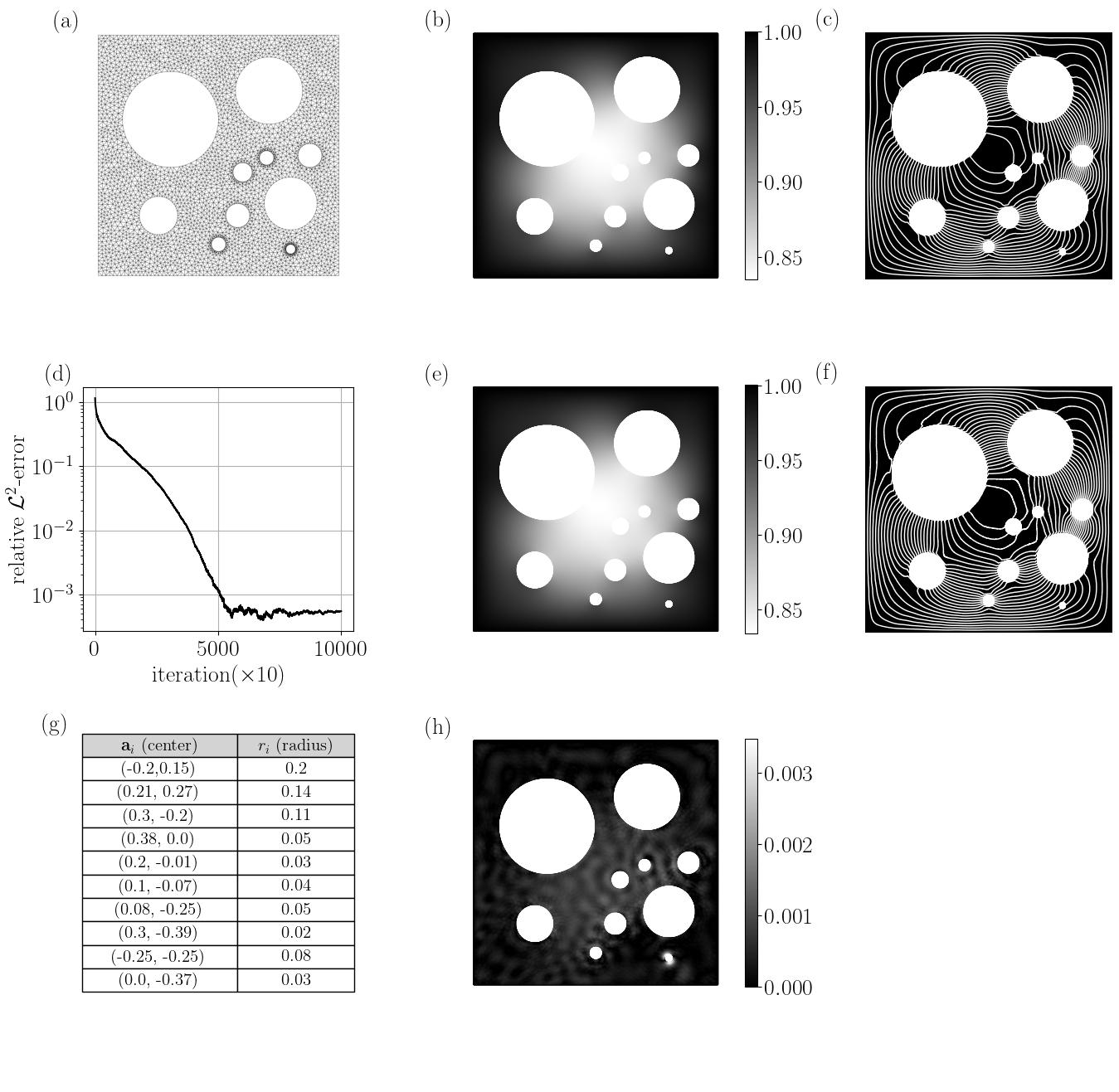}
\caption{Result of various perforation sizes. (a): Finite element mesh, (b): FEM solution, (c): contour map of the FEM solution, (d): learning trajectory of DFLM, (e): DFLM approximation, (f): contour map of the DFLM solution, (g): pointwise difference.}
\label{fig:example2}
\end{figure}

In the second example, we consider a collection of circular perforations of various sizes whose radii range from 0.02 to 0.2 (the location and size of the holes are detailed in Fig.~\ref{fig:example2}-(g)). The mesh around each circle becomes extremely small as the radius decreases (see Fig.~\ref{fig:example2}-(a) for the finite element mesh). We choose this configuration to check the performance of DFLM in handling various scales. 

As there is a broader range of scales than in the previous example, the Fourier embedding uses a more complex layer of dimension 400 and the variance $\sigma^2=25$ to account for broader scales. In contrast, the rest of the network structure remains the same: three hidden layers of dimension 200 with $\tanh$ activation function. Also, the perforation of radius 0.02 induces a small micro time step $\delta t=1\times 10^{-6}$ to estimate the exit time on the smallest perforation. That is, $\delta t$ is chosen to satisfy $\mathbb{E}[|\bm{B}_{\delta t}|]\ll 0.014<0.02$ (0.014 is related to a small perforation in the next example).

As the solution to this problem has more variations than the previous example, the macro time step, which is related to the trainability of the network \cite{DFLManalysis}, can be shorter than the previous case, which is set to $64\times \delta t$ (this value is ten times shorter than the previous case). Related to the empirical loss, there are $N_r=5000$ collocation points and $N_s=400$ stochastic walkers at each collocation point. Thus, the overall computational cost is comparable to the previous case. $N_r$ is 1.7 times more, but the number of steps used to sample the stochastic walkers is half the previous case (note that the number of stochastic walkers has not changed). The learning rate parameters are $(\alpha_0,\gamma)=(1.5\times 10^{-3},0.9)$

The first row of Fig.~\ref{fig:example2} corresponds to (a) the finite element mesh, (b) the reference FEM solution, and (c) the contour lines of the FEM solution (c). DFLM results in the second row match the reference results, including the contour lines. The total number of iterations is also comparable to the previous case; after about 60,000 iterations, the relative $\mathcal{L}^2$-error of DFLM stabilizes. The pointwise error is a bit larger in this case, which is of order $3\times 10^{-3}$ but the relative $\mathcal{L}^2$-error, $5.4557\times 10^{-4}$, is of the same order as in the previous case.

\subsection{Periodic perforations}
The last example is a homogenization problem with periodic perforations. Each perforation is of the same size, a circle of radius $1.428\times 10^{-2}$ regularly spaced in the square domain $\Omega$. The radius is chosen so that we can place 400 perforations $\{B(\bm{a}_{m,n};r): m,n=1,2,\cdots,20\}$ centered at a $20\times 20$ uniform grids $\bm{a}_{m,n}= \left(\frac{1}{40}+\frac{m}{20},\frac{1}{40}+\frac{n}{20}\right)$ (and thus leads to  $r= \frac{4}{7}\times \frac{1}{20}\times \frac{1}{2}\simeq 0.014$). In addition to the challenge imposed by the small radius, another challenge comes from the short average distance between perforations.

Resolving all relevant scales of the solution can be challenging. As DFLM aims at averaged behaviors of the solution, the solution structure can be simpler than the second example. Thus, we use the Fourier embedding of dimension 200 with variance $\sigma^2=9$ as in the first experiment, followed by three hidden layers of dimension 200. Another difference is the activation function, which is ReLU in this case. The number of collocation points is $N_r=5000$ while the number of stochastic walkers is $N_s=400$. Learning parameters are $(\alpha_0,\gamma)=(10^{-3},0.9)$.
As in the second experiment, the micro time step $\delta t$ is set to $1\times 10^{-6}$ to satisfy $\mathbb{E}[|\bm{B}_{\delta t}|] \ll 0.014$ and the macro time step $\Delta t$ is $64\times \delta t$. That is, the overall computational cost of this experiment is comparable to the two previous experiments. 

\begin{figure}[t!]
\centering
\includegraphics[width=1.0\textwidth]{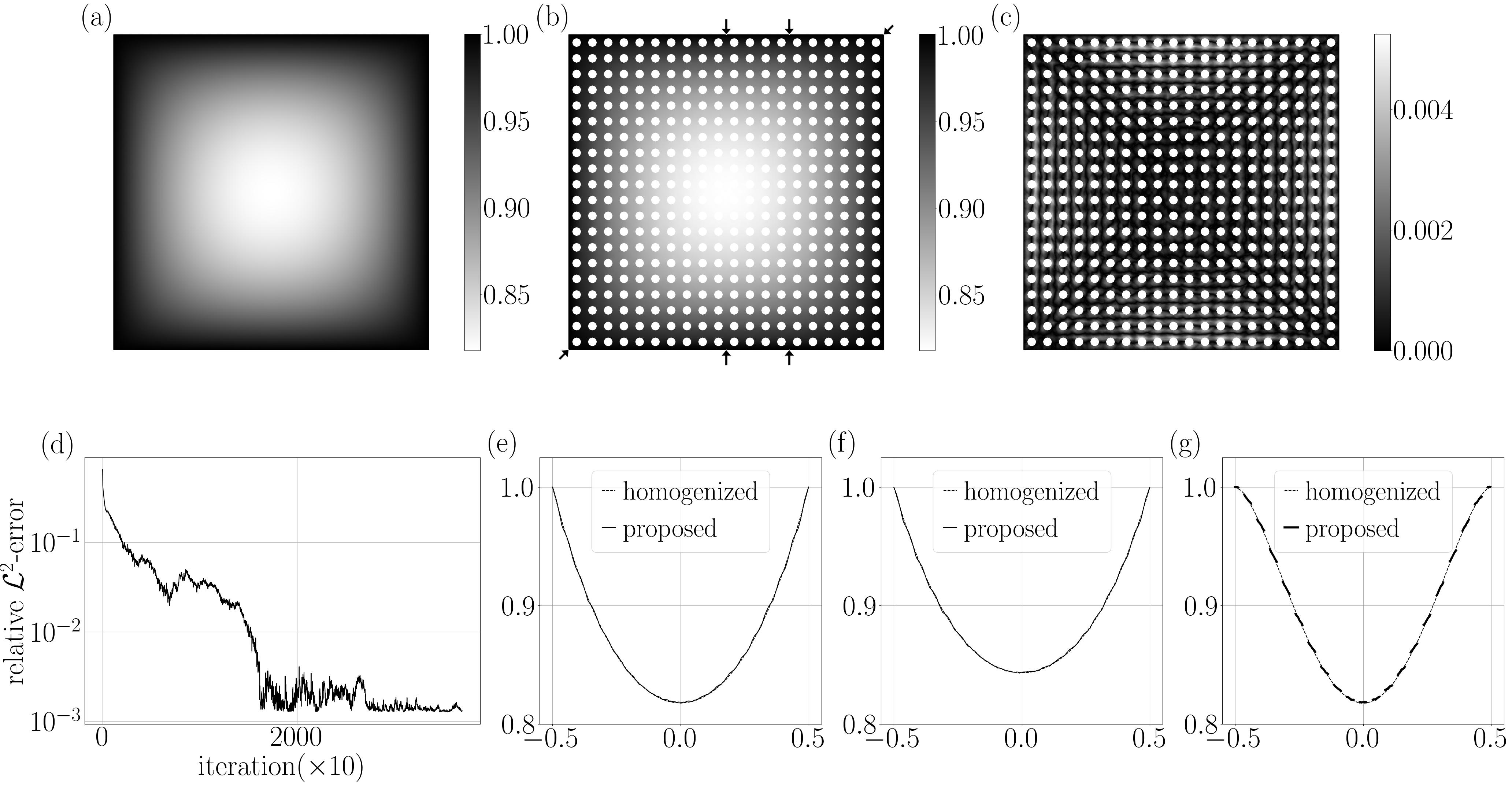}
\caption{Result of homogenization problem. (a): Homogenized solution, (b): DFLM solution, (c): pointwise error, (d): learning trajectory of DFLM, (e): cross-section at $x_1=0$, (f): cross-section at $x_2=0.2$, (g): diagonal cross-section at $x_1=x_2$. The arrows in (b) represent the directions of the cross-sections.}\label{fig:example3}
\end{figure}
In this stringent test case, it is challenging to obtain a FEM solution with convergence. Instead, we compare the DFLM result with the homogenized solution. Using the given perforation geometry in a periodic cell, we calculate the following homogenized tensor by numerically solving Eq.~\eqref{eq:homogenizedcoeff}
\begin{equation}
A^{0} \simeq 
\begin{bmatrix}
0.80783 & -1.56838\times 10^{-5} \\
-1.56838\times 10^{-5} & 0.80783
\end{bmatrix}.
\end{equation}
Using this homogenized tensor, we use FEM to calculate the homogenized solution (Fig.~\ref{fig:example3}-(a)). DFLM result (Fig.~\ref{fig:example3}-(b)) shows a global behavior of the solution matching the reference homogenized solution. Fig.~\ref{fig:example3}-(e)$\sim$(g) show the cross sections at $x_1=0$ (e), at $x_2=0.2$ (f), and at $x_1=x_2$ (g). The first two cross-sections bypass the perforations; thus, the DFLM result does not have missing parts. The DFLM results are on top of the reference solution in all three cross-sections with almost indistinguishable overlaps. The pointwise error is of order $5\times 10^{-3}$ with a relative $\mathcal{L}^2$-error $1.3034 \times 10^{-3}.$
We emphasize that DFLM non-intrusively captures the averaged behavior without knowing the homogenized tensor or calculating it through cell problems.

%
%
%
\section{Discussions and Conclusions}\label{sec:discussion}
We have extended the derivative-free loss method (DFLM) as a mesh-free approach for handling perforated domain problems, focusing on the efficient representation of multiscale solutions. As the boundary of the perforation is typically described as a homogeneous Neuman boundary condition, we have also implemented the Neumann boundary condition in DFLM using the reflecting stochastic process. As DFLM uses averaging to represent a PDE solution, a martingale representation (Eq.~\eqref{eq:q_martingale_stochastic_process} or Eq.~\eqref{eq:q_tilde_martingale_process}), it shows robust results in capturing macroscopic behaviors of multiscale solutions while maintaining the computational cost comparable to a non-multiscale solution. 

The current study has focused on an elliptic-type PDE in the form of Eq.~\eqref{eq:model}, which will provide a tool for investigating the thermal properties of perforated materials. Additionally, fluid motions through perforated materials have several applications in environmental science, such as air and water filtrations. DFLM has recently been extended to viscous fluid problems \cite{DFLMfluid}. We plan to investigate DFLM for fluid problems in perforated domains, considering various mathematical issues related to the time stepping of the micro and macro time steps that affect the trainability of the network. 

Also, we plan to consider perforation shapes other than circles. In such an extension, the issue is mainly in detecting the reflection point on the boundary. We believe this issue is technical but straightforward to investigate, which we leave as future work. As with other network-based approaches to solving PDEs, there are several unanswered questions we need to address, such as network complexity and design, along with several other tuning parameters. It would be natural to speculate a general guideline to answer these questions without much tuning procedure.

\section*{Acknowledgments}
The authors are supported by ONR MURI N00014-20-1-2595.

\bibliographystyle{siam}
\bibliography{dflm_perforated}
\end{document}